\documentclass[a4paper,12pt,english]{article}

\usepackage[english]{babel}
\usepackage[utf8]{inputenc}
\usepackage[T1]{fontenc}

\usepackage{float}
\usepackage{graphics}
\usepackage{graphicx}
\usepackage{epsfig}
\usepackage{indentfirst}
\usepackage{amsmath, amsfonts, amssymb, amsthm, mathtools}
\usepackage{url}
\usepackage{csquotes}
\usepackage{enumerate}

\usepackage[backend=biber, style=numeric-comp]{biblatex}
\DeclareTextFontCommand{\emph}{\boldmath\bfseries}
\DefineBibliographyStrings{english}{mathesis = {Master dissertation}}

\newcommand{\p}{\partial}
\newcommand{\E}{\mathrm{E}}

\begin{document}

\title{ANN-MoC Method for Solving Unidimensional Neutral Particle Transport Problems}

\author{Pedro H.A. Konzen\\
  \small{IME/UFRGS, Porto Alegre, RS, Brazil}\\
  \small{\texttt{pedro.konzen@ufrgs.br}}}

\date{\today}
\maketitle

\begin{abstract}
Neutral particle transport problems are fundamental in the modeling of energy transfer by radiation (photons) and by neutrons with many important applications. In this work, the novel ANN-MoC method for solving unidimensional neutral particle transport problems is presented. Following the Method of Discrete Ordinates (DOM) and decoupling with a Source Iteration (SI) scheme, the proposed method applies Artificial Neural Networks (ANNs) together with the Method of Characteristics (MoC) to solve the transport problem. Once the SI scheme converges, the method gives an ANN that estimates the average flux of particles at any points in the computational domain. Details of the proposed method are given and results for two test cases are discussed. The achieve results show the potential of this novel approach for solving neutral particle transport problems.

\end{abstract}

\section{Introduction}

Photon and neutron transport are important examples of neutral particles transport phenomena. The first appears in many applications, mainly in that involving energy transport via radiative transfer \cite{Modest2013a}. Practical applications includes the design of industrial furnaces, combustion chambers, or forming processes such as glass and ceramics manufacturing \cite{Frank2004a, Larsen2002a, Viskanta1987a}. Other applications are found in the fields of astrophysics \cite{Meinkohn2002a, Richling2001a}, medical optics \cite{Abdoulaev2003a,  Hielscher1998a, Tarvainen2008a, Wang2007a}, developing of micro-electro-mechanical systems \cite{Knackfuss2006a}. Neutron transport also has applications in medicine and clearly in nuclear energy generation \cite{Lewis1984a, Stacey2007a}. 

In this work, the neutral particle transport is assumed to be modeled in a unidimensional space domain $\overline{\mathcal{D}} = [a, b]$ as it follows
\begin{subequations}\label{eq:tp}
  \begin{align}
    &\forall \mu\in [-1, 1]:~\mu\cdot\frac{\p}{\p x} I(x,\mu) + \sigma_tI = \frac{\sigma_s}{2}\int_{-1}^1I(x,\mu')\,d\mu' + q(x,\mu), \forall x\in\mathcal{D},\label{eq:te}\\
    &\forall \mu>0:~I(a,\mu) = I_{a},\label{eq:tbca}\\
    &\forall \mu<0:~I(b,\mu) = I_{b},\label{eq:tbcb}
  \end{align}
\end{subequations}
where $I(x,\mu)$ is the angular flux of particles at the point $x\in\overline{\mathcal{D}}=[a, b]$ and in the direction $\mu\in [-1, 1]$, $\sigma_t$ is the total absorption coefficient and $\sigma_s$ the scattering coefficient, $q(x,\mu)$, $I_a$ and $I_b$ are, respectively, the sources in $\mathcal{D}$ and on its boundary. The average flux of particles is given by
\begin{equation}
  \Psi(x) := \frac{1}{2}\int_{-1}^1I(x,\mu)\,d\mu.
\end{equation}

Many solution approaches are available to problem \eqref{eq:tp} (see, for instance, \cite{Lewis1984a, Modest2013a}). One of the most applied is the so called Discrete Ordinates Method (DOM, \cite{Modest2013a}). By considering a numerical quadrature $\{\mu_i,w_i\}_{i=1}^{N}$, the problem \eqref{eq:tp} is approximated by a system of equations only for the discrete directions $\mu_i$, $i=1,2,\dotsc,N$. The equations can be further decoupled by using the Source Iteration (SI) strategy, where the system is iteratively solved for approximations of $\Psi(x)\approx \Psi^{(j)}(x)$, $j=1,2,3,\ldots$, until a given stop criteria. At each SI iterate, one has a decoupled system of $N$ linear first order partial differential equations, which can be solved by the Method of Characteristics (MoC, \cite{Evans2010a}). To do so, one will need to compute an integral depending on the $\Psi$ approximation.

In this work, we present a novel method to solve \eqref{eq:tp}, it integrates and Artificial Neural Network (ANN, \cite{Goodfellow2016a, Haykin2005a}) into the DOM-MoC approach. The main idea is to train an ANN to estimate the average flux $\Psi^{(j)}$ at each SI iterate. It is a meshless method, in the sense that it does not rely on a fixed domain mesh. After convergence, the method gives an ANN that estimate $\Psi(x)$ for all $x\in\overline{\mathcal{D}}$.

\section{The ANN-MoC Method}

Following the Discrete Ordinates Method (DOM), we assume a numerical quadrature $\{\mu_i, w_i\}_{i=1}^N$, and the Source Iteration (SI) approximation of problem \eqref{eq:tp} is given as follows
\begin{subequations}\label{eq:domp}
  \begin{align}
    &i=1,\dotsc,N:~\mu_i\cdot\frac{\p}{\p x} I^{(j)}(x,\mu_i) + \sigma_tI^{(j)}_i(x) = \sigma_s\Psi^{(j-1)}(x) + q(x,\mu_i), \forall x\in\mathcal{D},\label{eq:dome}\\
    &\mu_i>0:~I^{(j)}_i(a) = I_{a},\label{eq:dombca}\\
    &\mu_i<0:~I^{(j)}_i(b) = I_{b},\label{eq:dombcb}
  \end{align}
\end{subequations}
where $I^{(j)}_i \approx I^{(j)}(x,\mu_i)$, $l=1,2,\dotsc,L$, and $\Psi^{(0)}(x)$ is a given initial approximation for $\Psi(x)$. Then, the $j$-th approximation of the average flux is given by
\begin{equation}
  \Psi^{(j)}(x) = \frac{1}{2}\sum_{i=1}^N w_iI^{(j)}_i(x)
\end{equation}

Now we use the Method of Characteristics (MoC) by applying the change of variables $x(s) = x_0 + s\cdot\mu_i$. Then, for each $i=1,\dotsc,N$, equation \eqref{eq:dome} can be rewritten as follows
\begin{equation}
  \frac{d}{ds}I_i^{(j)}(s) + \sigma_t I^{(j)}_i(s) = \sigma_s\Psi^{(j-1)}(s) + q(s,\mu_i),
\end{equation}
where $I_i^{(j)}(s) = I_i^{(j)}\left(x(s)\right)$, and analogous for the other term. An integrating factor than gives us
\begin{equation}\label{eq:moc}
  I_i^{(j)}(s) = I_i^{(j)}(0)e^{-\int_0^s\sigma_t\,ds'} + \int_0^s\left[\Psi^{(j)}(s')+q(s',\mu_i)\right]e^{-\int_{s'}^s\sigma_t\,ds''}\,ds'
\end{equation}
The computation of the integral term involving $\Psi^{(j)}(s)$ is an issue, since it usually requires the evaluation of $\Psi^{(j)}(s)$ at several points $s\in (0, s)$, which can be a large interval depending on the direction $\mu_i$.

The idea of the proposed ANN-MoC method, is to train an Artificial Neural Network (ANN) to estimate $\Psi^{(j)}$ at each source iteration. In the following, we simplify the notation by omitting the super-index $(j)$.

\subsection{ANN Average Flux Estimation}

The ANN is assumed to be a Multilayer Perceptron (MLP, \cite{Haykin2005a}) that has $x\in\overline{\mathcal{D}}$ as input and the estimate $\tilde{\Psi}(x)$ as output. It is denoted by
\begin{equation}\label{eq:ann}
  \tilde{\Psi}(x) = \mathcal{N}\left(x; \{(W^{(l)},\pmb{b}^{(l)}, f^{(l)})\}_{l=1}^{n_l}\right),
\end{equation}
where $(W^{(l)},\pmb{b}^{(l)}, f^{(l)})$ denotes the triple of the weights $W^{(l)} = \left[w^{(l)}_{i,j}\right]_{i,j=1}^{n^{(l-1)}, n^{(l)}}$, the bias $\pmb{b}^{(l)} = \left(b^{(l)}_i\right)_{i=1}^{n^{(l)}}$ and the activation function $f^{(l)}$ in the $l$-th layer of the network. The number of neurons (units) at each layer is denoted by $n^{(l)}$, $l=1,2,\dotsc,n_l$. The MLP forwardly computes
\begin{align}
  \pmb{a}^{(l)} = f^{(l)}\left(W^{(l)}\pmb{a}^{(l-1)}+\pmb{b}^{(l)}\right),
\end{align}
where $\pmb{a}^{(0)} = x$ and $\tilde{\Psi}(x) = \pmb{a}^{(n_l)}$.

Given a fixed structure (number of layers $n_l$, number of units $n^{(l)}$ per layer and the activation functions), the training of the ANN consists in solving the following optimization problem
\begin{equation}\label{eq:annoptim}
  \min_{\{(W^{(l)},\pmb{b}^{(l)})\}_{l=1}^{n_l}}\frac{1}{n_s}\sum_{m=1}^{n_s}\left(\tilde{\Psi}^{(m)}-\Psi^{(m)}\right)^2
\end{equation}
for a given training set $\{x^{(m)}, \tilde{\Psi}\left(x^{(m)}\right)\}_{m=1}^{n_s}$, where $n_s$ is the number of samples.

\subsection{The ANN-MoC Algorithm}

The proposed ANN-MoC method computes successive approximations of the average flux $\Psi(x)$ for all points in the domain $\overline{\mathcal{D}}$. It starts from the ANN \eqref{eq:ann} trained with given initial training set $\{x^{(m)}, \tilde{\Psi}^{(0)}\left(x^{(m)}\right)\}_{m=1}^{n_s}$, for randomly selected points $x^{(m)}\in\overline{\mathcal{D}}$, $m=1,2,\dotsc,n_s$. Then, the approximation $\tilde{\Psi}^{(j)}$ is iteratively computed from its previous $\tilde{\Psi}^{(j-1)}$ by solving the problem \eqref{eq:domp} from the MoC solution \eqref{eq:moc} and by replacing $\Psi^{(j)}(s')$ for its estimate from the ANN $\mathcal{N}(s')$, trained on the last $l-1$-th source iteration.

The ANN-MoC algorithm follows the steps:
\begin{enumerate}[1.]
\item Set the ANN structure $\mathcal{N}(x)$ with random weights and bias.
\item Set an initial approximation $\Psi^{(0)}(x)$ for all $x\in\overline{\mathcal{D}}$.
\item Set $n_s$ and the set of points $\{x^{(m)}\}_{m=1}^{n_s}$.
\item Train the ANN with the training set $\{x^{(m)}, \Psi^{(0)}\left(x^{(m)}\right)\}_{m=1}^{n_s}$.
\item Set the quadrature $\{\mu_i, w_i\}_{i=1}^{N}$.
\item For $j=1,\dotsc,L$:
  \begin{enumerate}[6.a)]
  \item For $i=1,\dotsc,N$, for $m=1,\dotsc,n_s$:
    \begin{itemize}
    \item If $\mu_i>0$, then $s=(x^{(m)}-a)/\mu_i$
      \begin{equation}\label{eq:mocp}
        I^{(j)}_i(x^{(m)}) = I_ae^{-\int_0^s\sigma_t\,ds'} + \int_0^s\left[\mathcal{N}(s')+q(s',\mu_i)\right]e^{-\int_{s'}^s\sigma_t\,ds''}\,ds'
      \end{equation}
    \item If $\mu_i<0$, then $s=(x^{(m)}-b)/\mu_i$
      \begin{equation}\label{eq:mocn}
        I^{(j)}_i(x^{(m)}) = I_be^{-\int_0^s\sigma_t\,ds'} + \int_0^s\left[\mathcal{N}(s')+q(s',\mu_i)\right]e^{-\int_{s'}^s\sigma_t\,ds''}\,ds'
      \end{equation}
    \end{itemize}
  \item Compute $\Psi^{(j)} = \frac{1}{2}\sum w_iI_i^{(j)}$.
  \item Retrain the ANN $\mathcal{N}(x)$ with the new training set $\{x^{(m)}, \Psi^{(j)}\left(x^{(m)}\right)\}_{m=1}^{n_s}$.
  \item Check a given stop criteria.
  \item Reset the random set of points $\{x^{(m)}\}_{m=1}^{n_s}$.
  \end{enumerate}
\end{enumerate}

\section{Results}

In this section we present results of the application of the ANN-MoC method to solve two different problems. The first is set from a manufactured solution and the second is a benchmark problem selected from the specialized literature.

\subsection{Problem 1: Manufactured Solution}

We assume the exact angular fluxes are given as
\begin{equation}
  \hat{I}(x,\mu) = e^{-\alpha\sigma_t x}.
\end{equation}
By substituting in \eqref{eq:te}, one obtains the source
\begin{equation}
  q(x,\mu) = \left(\kappa - \alpha\sigma_t\mu)\right)e^{-\alpha\sigma_t x}.
\end{equation}
The exact average particle flux can be also analytically calculated as
\begin{equation}
  \hat{\Psi}(x) = e^{-\alpha\sigma_t x}.
\end{equation}

\begin{figure}[H]
\centering
\includegraphics[width=.7\textwidth]{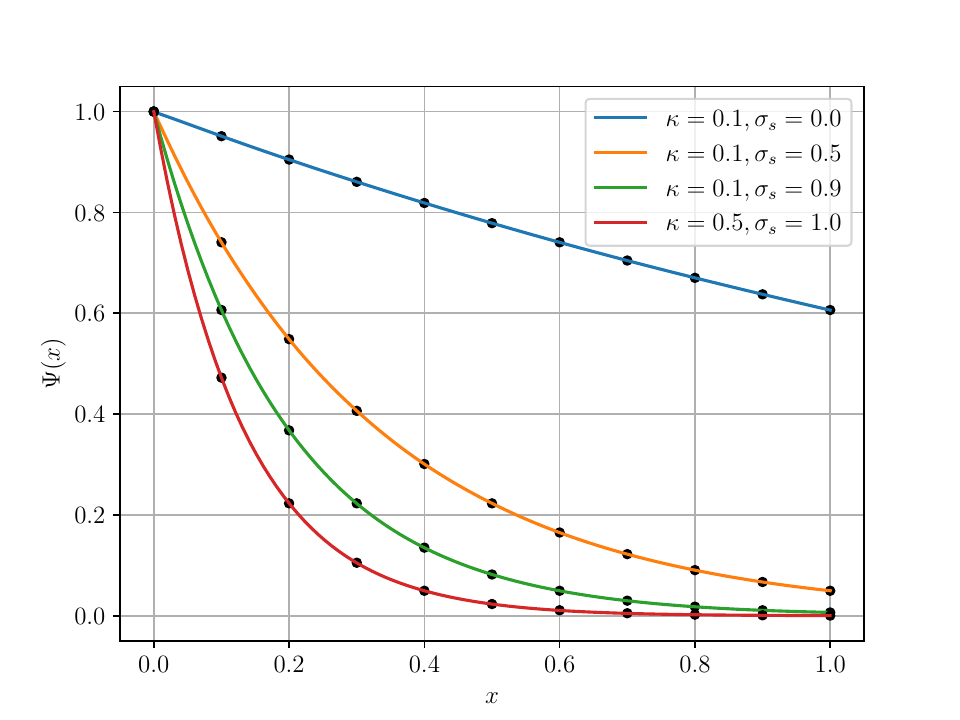}
\caption{{\small ANN-MoC (lines) {\it versus} the exact (dots) solutions of Problem 1 for several values of $\kappa$ and $\sigma_s$.}}
\label{fig:ms}
\end{figure}

Figure \ref{fig:ms} shows a comparison of the ANN-MoC {\it versus} the exact solutions for Problem 1 with several different values of $\kappa$ and $\sigma_s$. The approximated solutions have been achieved by using an 1-100-50-5-1 MLP with hyperbolic tangent as activation function on hidden layers and the sigmoid function to activate the output neuron. The Adam method \cite{Kingma2017a} has been used for solving the optimization problem \eqref{eq:annoptim} at each training step. The Gauss-Legendre quadrature with $N=100$ nodes has been assumed for the DOM angular discretization and the number of point samples has been fixed to $n_s=101$. As stop criteria for the SI iterations, we have applied
\begin{equation}
  \|\tilde{\Psi}^{(l)}-\tilde{\Psi}^{(l-1)}\|_2 < \max\{\varepsilon, \varepsilon\|\tilde{\Psi}^{(l)}\|\},
\end{equation}
with $\varepsilon = 10^{-5}$.

\begin{table}[H]
  \caption{{\small Average flux of particles computed at selected points for Problem 1 with $\kappa=0.1$ and $\sigma_s=0.5$.}}
  \centering
  \begin{tabular}{r|ccccc|r}\hline
    $n_s$ & $\Psi(0.0)$ & $\Psi(0.25)$ & $\Psi(0.5)$ & $\Psi(0.75)$ & $\Psi(1.0)$ & $\|\tilde{\Psi}-\hat{\Psi}\|_2$\\\hline
    11    & 1.0000 & 0.4722 & 0.2232 & 0.1051 & 0.0498 & $1.98\E-4$\\
    51    & 0.9998 & 0.4726 & 0.2232 & 0.1053 & 0.0498 & $1.48\E-4$\\
    101   & 0.9992 & 0.4724 & 0.2231 & 0.1053 & 0.0496 & $1.07\E-4$\\
    201   & 0.9995 & 0.4722 & 0.2231 & 0.1054 & 0.0496 & $1.22\E-4$\\\hline
    exact & 1.0000 & 0.4724 & 0.2231 & 0.1054 & 0.0498 & -x-\\\hline 
  \end{tabular}
  \label{tab:prob1}
\end{table}

Table \ref{tab:prob1} presents the average flux of particles computed at selected domain points for Problem 1 with $\kappa=0.1$ and $\sigma_t=0.5$. One can observe that the increase of sample points from $n_s=11$ to $201$ produce similar results, which indicates the training of the MLP will not profit from further increasing the number of samples. This is due to the randomization of the sample points at each SI iteration.

\subsection{Problem 2: Benchmark Solution}

The second application of the ANN-MoC is for the benchmark problem available in the work \cite[Table 1]{Vargas2017a}. The problem sources are
\begin{equation}
  q(x,\mu) = x - x^2,
\end{equation}
and $I_a=I_b = 0$. The absorption coefficient is fixed to $\sigma_t=1$.

\begin{figure}[H]
\centering
\includegraphics[width=.7\textwidth]{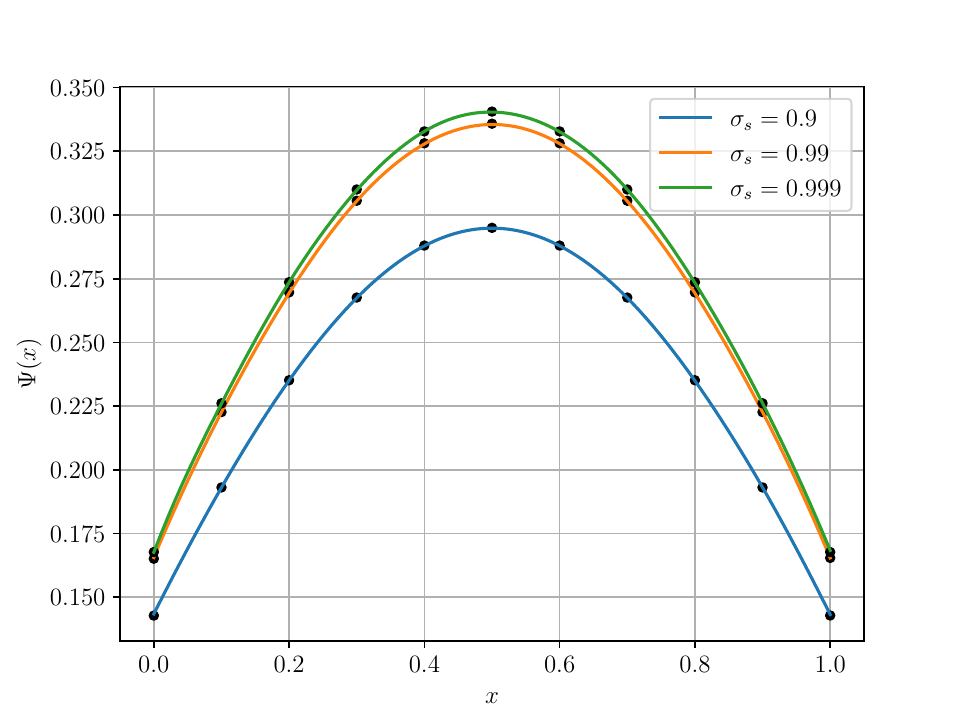}
\caption{{\small ANN-MoC (lines) {\it versus} the exact (dots) solutions for Problem 2 with $\sigma_s=0.9, 0.99$ and $0.999$.}}
\label{fig:bs}
\end{figure}

Figure \ref{fig:bs} shows a comparison of the ANN-MoC (lines) {\it versus} the exact (dots) solutions for Problem 2 with the scattering coefficient set to $\sigma_s=0.9, 0.99$ and $0.999$. The ANN-MoC parameters were all set as the same used for solving the Problem 1, with $n_s=101$. As in that case, we can observe very good accordance between the proposed method and the expected solutions.

\section{Final Considerations}

In this paper, the novel ANN-MoC method has been presented for solving unidimensional neutral particle transport problems. Its main idea is to apply an ANN for the estimates of the average flux of particles computed from a DOMM-MoC approach. One of its advantages is to be a meshless method, since no fixed mesh is necessary in the computations. After the convergence of the SI iterations, the method gives an ANN to estimate the average flux at any point of the domain. The achieved first results have been presented and they show a very good accordance between the ANN-MoC and the expected solutions. This indicates the potential of the method as an alternative to be applied for the solution of more complex transport problems. Further work should also address on the ANN-MoC comparison with the classical strategy of estimating the average fluxes by interpolation on mesh points. If from one point of view the training and evaluation of an ANN is more expensive to compute that performing interpolation, it may be compensated by the need of relatively small number of sample points on a meshless structure.

\printbibliography

\end{document}